\documentclass[11pt,a4paper]{report}
\usepackage{amsmath}
\usepackage{setspace}
\usepackage{xypic}
\usepackage[all]{xy}
\usepackage{latexsym}
\usepackage{theorem}
\newtheorem{teorema}{Theorem}[section]
\newtheorem{definicion}[teorema]{Definition}
\include{amslatex}
\newtheorem{proposicion}[teorema]{Proposition}

\newtheorem{corolario}[teorema]{Corollary}
{\theorembodyfont{\rmfamily}
}
{\theorembodyfont{\rmfamily}
}

\numberwithin{equation}{section}

\include{amsmath}

\begin{document}
\begin{title}
{\LARGE {\bf Some Rigidity Conditions on Berwald Structures}}
\end{title}
\maketitle
\author{

\begin{center}
Ricardo Gallego Torrome,\\
\paragraph{}
\paragraph{}
\paragraph{}
Department of Physics, Lancaster University,
\\
Lancaster, LA1 4YB, United Kingdom \& The Cockcroft Institute, UK\\[3pt]

\paragraph{}

Departamento de Matem{a}ticas, Estad{\'i}stica y Computaci\'{o}n,
\\
Facultad de Ciencias, Universidad de
Cantabria, \\

Avda. de los Castros, s/n, 39071 Santander, Spain

\end{center}}

\tableofcontents{}

\newpage
\paragraph{}
\paragraph{}
\paragraph{}
This work is dedicated to my parents, Luis y Elvira and to my
brothers Elvira Julio and Luis Miguel.
\newpage
\paragraph{}
\paragraph{}
\paragraph{}
The present manuscript corresponds to the memory of my research work in order to obtain the Masters Degree
in Mathematics. It has been written under the supervision
of Prof. Fernando Etayo Gordejuela. I am pleased to thank to him
for his time, his kind hospitality, his
collaboration during the elaboration of the present and related works. Undoubtable, this thesis can not be done
without his guidance. Therefore, my gratitude to him.
\newpage
\paragraph{}
\begin{center}
{\bf Sobre una condicion de rigidez de los espacios de Berwald}
\end{center}
\paragraph{}
\paragraph{}
{\bf Resumen}. Esta tesis contiene una introduci\'{o}n al
m\'{e}todo de los promediados de estructuras geometricas, en
particular de estructuras definidas por conexiones Finslerianas.
Se aplica el m\'{e}todo a espacios de Berwald, que son espacios de
Finsler pero que preservan todavia mucha de las
caracter\'{i}sticas propias de los espacios Riemannianos. En este
sentido, se obtienen condiciones de rigidez geod\'{e}sica, como el
{\it teorema} 5.1.3. En la prueba, es esencial el promediado de la
conexi\'{o}n the Chern. Mas tarde se muestra que la conexi\'{o}n
de Levi Civita de cualquier m\'{e}trica Riemanniana af\'{\i}nmente
equivalente a una estructura de Berwald deja invariante por
transporte paralelo la indicatriz de dicha estructura de Berwald.
Tambi\'{e}n se demuestra el resultado rec\'{i}proco: Si $({\bf
M},F)$ es una estructura de Finsler y existe una estructura
Riemanniana cuya conexi\'{o}n de Levi Civita deja invariante por
transporte paralelo la indicatriz de la estructura de Finsler,
entonces $({\bf M},F)$ es de Berwald. Como aplicaci\'{o}n se
obtiene una condici\'{o}n necesaria para que una variedad sea de
Landsberg pura.

\chapter{Introduction}

Finsler geometry has its conceptual genesis in the seminar {\it
lecture} or Berhand Riemann "Ueber die Hypothesen, welche der
Geometrie zu Grunde liegen" (Habilitationsschrift, 1854,
Abhandlungen der Königlichen Gesellschaft der Wissenschaften zu
Göttingen, 13 (1868)). In this work, Riemann introduced the basic
ingredients of the modern notion of manifold, Riemannian structure
and Finsler structure. However, in the same work he noticed the
complication of the (general) Finsler case compared with the
Riemannian (quadratic case).

Due to the pre-eminence of the quadratic case, Finsler Geometry
was dormmitant for decades, reappearing in the thesis of P.
Finsler (under the supervision of Caratheodory) in 1918. This is
the reason why this type of geometry is also known by
$Riemann-Finsler$ geometry or Finsler Geometry for short.

After Finsler's thesis, a explosion in the field came in the next
decades and diverse schools of Finsler Geometry emerged, as well as
significant contribution of many geometers. During this earlier
development, the results are mainly local in character and related
with analytical questions, in particular, the calculus of
variations.

One of the relevant figures working on that period on Finsler
Geometry was L. Berwald, who introduced a connection and a class
of spaces sharing his name. Berwald connection is important
because can be extracted directly form the differential equation
stipulated as being the geodesic equations. Berwald spaces are
interesting because they are closely related to Riemannian spaces.
Berwald spaces is the category that, because being quite close to
the Riemannian category, more rigidity conditions can be found.
Indeed, a rigidity result due to Szab\'o says that a Berwald space
of dimension $2$ is either a Riemannian space or a locally
Minkowski space ([6]); therefore to find examples of Berwald
spaces it is necessary to go higher dimensions ([1, {\it chapter
11}]). In addition, Berwald spaces have the interesting feature
that they are related to the {\it Equivalence Principle} of General Relativity:
 the Berwald connection of a Berwald spaces constitute a
general type of torsion-free connections compatible with it.

It is in the category of Berwald Spaces where the present thesis
has to be considered. This memory explains in a (hopefully)
self-contained way some of the results presented in reference [15]
in a jointly work with Prof. Fernando Etayo.

Between the amount of results presented, we would like to mention
the following:
\begin{enumerate}
\item A result on geodesic rigidity in the category of Berwald
spaces ({\it proposition 5.1.7}), which is similar to a result
obtained by V. Matveev ([9])

\item A rigidity result on Berwald spaces, ({\it proposition
5.2.4}).

\item A rigidity condition for Landsberg spaces {\it proposition
5.3.2}.
\end{enumerate}

The technical tool used to obtain these results is to consider the
average of some geometric Finslerian quantities and in particular,
the average of some Finsler linear connection and the averaged of
the fundamental tensor. The averaging operation is presented in
the pre-print [2] as well as in sub-sequent works. We think that
the averaged founded in that reference is useful (for some
purposes more than the averaged of the fundamental tensor) because the
relation between the average of the curvature of the original
connection and the curvature of the averaged connection. This
makes this average more powerful than the average of the
fundamental tensor.

\chapter{Basic Notions on Riemann-Finsler Geometry}

Let $(x,{\bf  U})$ be a local coordinate system on ${\bf M}$, where $x\in 
{\bf U}$ have local coordinates $(x^{1},...,x^{n})$, ${\bf U}\subset {\bf M}$ is an open set 
and ${\bf TM}$ is the tangent bundle of the manifold ${\bf M}$. A
tangent vector at the point $x\in {\bf M}$ is denoted by
$y^{i}\frac{{\partial}}{{\partial} x^{i}}\in {\bf T}_x {\bf M},\, y^i \in {\bf R}$.
We use Einstein's convention 
for up and down equal indices in this work if the contrary is not
directly stated.  We
can identify the point $x$ with its coordinates $(x^{1},...,x^{n})$ and the tangent vector $y\in {\bf T}_x 
{\bf M}$ at $x$ with its components $y=(y^1 ,...,y^n )$. Then each
local coordinate system $(x,{\bf U})$ on the manifold {\bf M}
induces a local coordinate system on ${\bf TM}$ denoted by
$(x,y,{\bf U})$ such that $y=y^{i}\frac{{\partial}}{{\partial}
x^{i}}\in {\bf T}_x{\bf M}$ has local natural coordinates
$(x^1,...,x^n ,y^1,...,y^n)$ in the induced natural coordinate
system. The slit tangent bundle is  $\pi:{\bf N}\longrightarrow
{\bf M}$ such that ${\bf N}={\bf TM}\setminus \{0\}$; i.e., the
tangent bundle with the zero section removed.
\section{Definition of Finsler Structure}

\begin{definicion}

A Finsler structure $F$ on the manifold ${\bf M}$ is a
non-negative, real function  $F:{\bf TM}\rightarrow [0,\infty [$
such that
\begin{enumerate}
\item It is smooth in the slit tangent bundle ${\bf N}$.

 \item Positive homogeneity holds: $F(x,{\lambda}y)=\lambda F(x,y)$ for every $\lambda >0$.

 \item Strong convexity holds:
the Hessian matrix
 \begin{equation}
g_{ij}(x,y): =\frac{1}{2}\frac{{\partial}^2 F^2
(x,y)}{{\partial}y^i {\partial}y^j }
\end{equation}
is positive definite on ${\bf N}$.
\end{enumerate}
\end{definicion}
We also denote by a Finsler structure to the pair $({\bf M},F)$. The coefficients $g_{ij}(x,y)$ are the 
components of the fundamental tensor $g$ defined later.

{\bf Remark 1}. When the second Bianchi identities are used,
the minimal 
smoothness requirement for the Finsler structure is to be $\mathcal{C}^5$; more generally, only 
$\mathcal{C}^4$ differentiable structure is required, if one speaks only of curvatures. 

{\bf Remark 2}. The homogeneity condition can be stronger: 
$F(x,\lambda y)=|\lambda |F(x,y)$. In this case $({\bf M},F)$ is called absolutely 
homogeneous Finsler structure.

{\bf Remark 3}. In some examples it is convenient to reduce the
condition of strong convexity in the whole ${\bf
TM}\setminus\{0\}$ to some proper sub-manifold of {\bf N} defined by
proper subsets of the tangent spaces. Then one speaks of
$y$-locality of the strong convexity condition.

\begin{definicion}$([1])$
Let $({\bf M}, F)$ be a Finsler structure and $(x,y,{\bf U}) $ a
local coordinate system induced on ${\bf TM}$ from the coordinate system $(x,{\bf U})$ of 
{\bf M}. The components of the Cartan tensor are defined by the
set of functions:
\begin{equation}
{ A}_{ijk}=\frac{F}{2} \frac{\partial g_{ij}}{\partial
y^{k}},\quad i,j,k=1,...,n.
\end{equation}
\end{definicion}
These coefficients are homogeneous of degree zero in $(y^1
,...,y^n )$. In the
Riemannian case the coefficients ${ A}_{ijk}$ are zero, and this fact
characterizes Riemannian geometry from 
other types of Finsler geometries (this result is known as
Deicke's theorem ([1, pg 393]).
\begin{center}
***
\end{center}

Let us consider the vector bundle $ {\pi}^* {\bf TM} $, pull-back
bundle of ${\bf TM}$ by the projection ${\pi}$, defined as the
minimal sub-bundle of the cartesian product ${\bf N}\times{\bf
TM}$ such that the following sub-bundle commutes:
\begin{displaymath}
\xymatrix{\pi^*{\bf TM} \ar[d]_{\pi_1} \ar[r]^{\pi_2} &
{\bf TM} \ar[d]^{\pi}\\
{\bf N} \ar[r]^{\pi} & {\bf M}}
\end{displaymath}
where the projection $\pi$ is
\begin{displaymath}
{\pi}:{\bf N}\longrightarrow {\bf M}
\end{displaymath}
\begin{displaymath}
u \longrightarrow x ,\, u\in{\bf N},\, x\in {\bf M}.
\end{displaymath}
${\pi}^* {\bf TM}$ has base manifold ${\bf N}$, the
fiber over the point $u\in {\bf N}$ with coordinates $(x,y)$ is isomorphic to ${\bf 
T}_{x}{\bf M}$ for every point $u\in \pi^{-1}(x)$ and the structure group is isomorphic to ${\bf GL}(n,{\bf 
R})$. Given a vector field $Z\in \Gamma {\bf TM}$, the
corresponding element on the pull-back bundle is defined on each
$u\in \pi^{-1}(x)\subset {\bf N}$ by the cartesian pair
$(u,Z(x))$.
\paragraph{}
An alternative treatment of Finsler geometry uses the homogeneity
properties on $y$ of the different geometric objects that appear
in the theory. In fact, for positive homogeneous metrics, one can
investigate the geometry of analogous pull-back bundles but where
the base manifold is the sphere bundle ${\bf SM}$ (or the
projective sphere bundle ${\bf PTM}$ in the case of absolutely
homogeneous structures). The sphere bundle ${\bf SM}$ is defined as
follows. Consider the manifold ${\bf N}$ and the equivalence
relation defined as
\begin{displaymath}
(x,y)\equiv (x,\tilde{y})\,\,\textrm{iff}\,\, \exists
\lambda\in {\bf R}^+\textrm{such that}y=\lambda \tilde{y}.
\end{displaymath}

 Then ${\bf SM}$ is a fiber bundle over the manifold 
${\bf M}$, with fiber over the point $x\in{\bf M}$
\begin{displaymath}
\pi: {\bf SM }\longrightarrow {\bf M}, \quad
(x,[y])\longrightarrow x,
\end{displaymath}
where $(x,[y])$ is the equivalence class defined by above
equivalence relation. Then, one can construct as before the
pull-back bundle ${\pi^* \bf TM}$. If the structure $F$
is absolutely homogeneous of degree zero, then one can define the
projective bundle in a similar way:
\begin{displaymath}
\pi _S : {\bf N}\longrightarrow {\bf PTM},
 \quad
(x,y)\longrightarrow (x,[y]),
\end{displaymath}
defining the equivalence class as $ [y ]:=\{(x,y)\,|\, y=\lambda
y_o,\,\,\forall \lambda\neq 0\,\}$.

For example, the matrix coefficients $ \big(g_{ij}(x,y)\big) $
are also invariant under a positive 
scaling of $y$ and therefore they live on ${\bf SM}$. The Cartan
tensor components ${ A}_{ijk}$ also live on ${\bf SM}$, if they
are defined according to the formula $(2.1.2)$. If $F$ is absolutely
homogeneous rather than positive homogeneous the coefficients
$g_{ij}$ and $A_{ijk}$ live on
${\bf PTM}$, the projective tangent bundle. ${\bf PTM}$ is defined in a similar way as ${\bf 
SM}$ but the projection also sends $y$ and $-y$ to the same
equivalence class $[y_0 ]$.
\begin{center}
***
\end{center}

{\bf Examples of Finsler Structures}
\begin{enumerate}
\item {\bf Minkowski Space}. Given a vector space ${\bf V}$ a Minkowski
norm is a map $\|,\|:{\bf V}\longrightarrow {\bf R}$, such that
\begin{enumerate}
\item It is non-negative and $\|y\|=0$ iff  $y=0$.

\item It is positive homogeneous of degree 1.

\item It is smooth on $y$ and the Hessian respect to $y$ is
strictly positive.
\end{enumerate}
A Minkowski space is a pair $({\bf V}, F)$ as above. Indeed, one can check that an ordinary norm
is also defined from the axioms of Minkowski norm.

\item {\bf Riemannian Structures.} In this case $F$ has the form
$F(x,y)=\sqrt{g_{ij}(x)y^iy^j},$, for each $x\in {\bf M}$ and
$y\in {\bf T}_x{\bf M}$ and the matrix $g_{ij}(x)$ defines  a positive definite, symmetric
bilinear form on ${\bf T}_x {\bf M}$.

\item $(\alpha, \beta)$-{\bf metrics}. They are Finsler structures
determined by a Riemannian norm $\alpha:=\sqrt{a_{ij}(x)y^i\,y^j}$
and a linear form $\beta:=\beta_i (x)y^i$. One of the most
interesting cases are Randers structure ([1,{\it chapter 11}]),
which has the form $F(x,y)=\alpha(x,y)+\beta(x,y)$. The $1$-form
$\beta$ has norm less than $1$ by the Riemannian norm $\alpha$.
This ensures positivity as well as strong convexity.

\item {\bf Numata} metrics. They are defined by functions of the
form $F(x,y)=\alpha(x,y)+\beta(x,y)$, where
$\alpha=\sqrt{g_{ij}(y)y^i y^j}$ is a homogeneous function of degree $1$.

\item It was proved that the function measuring the time spent
climbing a mountain can be represented by a Finsler function. One nice reference is [16]. Let us
consider the Finslerian distance between two arbitrary points $p$ and $d$, which is the infimun of the Finslerian length of all possible piecewise smooth path connecting them:
\begin{equation}
d(p,q):=\, \inf\,\{\int_{\sigma}ds \,\sqrt{\eta_{ij}\dot{\sigma}^i\dot{\sigma}^j},\,\sigma:[0,1]\longrightarrow {\bf M}\},
\end{equation}
where {\bf M } is a smooth representation of the mountain, $\eta$ is the inner Riemannian metric on {\bf M} induced from the Euclidean metric in
${\bf R}^3$, $\sigma$ is a path connecting $p$ and $q$ and $\dot{\sigma}$ the tangent vector along $\sigma$. Given a point over a possible path $\sigma(s)$, let us denote the maximal speed as $c(\sigma(s),\dot{\sigma}(s))$. Then, the minimal time is given by:
\begin{equation}
t_{min}(p,q):=\, \inf\,\{\int_{\sigma}ds \, \sqrt{\frac{\eta_{ij}\dot{\sigma}^i\dot{\sigma}^j}{c^2(\sigma(s),\dot{\sigma}(s))}},\,\sigma:[0,1]\longrightarrow {\bf M}\}.
\end{equation}
Since by definition $c(\sigma(s),\dot{\sigma}(s))$ is homogeneous of degree zero on the second argument, the function
\begin{displaymath}
F:{\bf N}\longrightarrow {\bf R}
\end{displaymath}
\begin{equation}
(x,y)\longrightarrow \sqrt{\frac{\eta_{ij}y^iy^j}{c^2(x,y)}}.
\end{equation}
This is a Finsler metric.
\end{enumerate}

\section{The Non-Linear connection}

An Ehresmann connection in a principal fiber bundle $\pi:{\bf
P}\longrightarrow {\bf M}$ is a splitting of the tangent bundle
${\bf TP}$ such that ${\bf T}_u {\bf P}=\mathcal{V}_u \oplus
\mathcal{H}_u$ with $\mathcal{V}_u =ker\, d\pi$, for all $u\in {\bf
P}$.

There is a non-linear connection on the manifold ${\bf N}$. In order
 to introduce it, let us define the {\it non-linear connection coefficients}, defined by the formula in local coordinates
\begin{displaymath}
\frac{N^{i}_{j}}{F}={\gamma}^{i}_{jk}\frac{y^{k}}{F}-A^{i}_{jk}
{\gamma}^{k}_{rs}\frac{y^{r}}{F}\frac{y^{s}}{F},\quad
i,j,k,r,s=1,...,n
\end{displaymath}
where the formal second kind Christoffel symbols
${\gamma}^{i}_{jk}$ are defined by the expression
\begin{displaymath}
 {\gamma}^{i}_{jk}=\frac{1}{2}g^{is}(\frac{\partial g_{sj}}{\partial
x^{k}}-\frac{\partial g_{jk}}{\partial x^{s}}+\frac{\partial
g_{sk}}{\partial x^{j}}),\quad i,j,k=1,...,n;
\end{displaymath}
$A^i _{jk}=g^{il}A_{ljk}$ and $g^{il}g_{lj}=\delta ^i _j .$ Note that the coefficients 
$\frac{N^{i}_{j}}{F}$ are invariant under positive scaling $y\rightarrow \lambda y$, $\lambda \in 
{\bf R^+ }$, $y\in {\bf T}_x {\bf M}$.

Let us consider the local coordinate system $(x,y,{\bf U}) $ of the manifold ${\bf TM}$. An induced 
tangent basis for ${\bf T}_u {\bf N},\, u\in {\bf N}$ is defined
by the vectors([2]):
\begin{displaymath}
\{ \frac{{\delta}}{{\delta} x^{1}}|_u ,...,\frac{{\delta}}{{\delta} x^{n}} |_u, 
F\frac{\partial}{\partial y^{1}} |_u,...,F\frac{\partial}{\partial
y^{n}} |_u\},
\end{displaymath}
\begin{equation}
 \frac{{\delta}}{{\delta} x^{j}}|_u =\frac{\partial}{\partial
x^{j}}|_u -N^{i}_{j}\frac{\partial}{\partial y^{i}}|_u ,\quad
i,j=1,...,n.
\end{equation}
The local sections $\{ \frac{{\delta}}{{\delta} x^{1}}|_u 
,...,\frac{{\delta}}{{\delta} x^{n}}|_u,\, u\in \pi^{-1}(x),\, x\in {\bf U} \} $ generates 
the local horizontal distribution $\mathcal{H}_U $, while $\{ \frac{\partial}{\partial 
y^{1}}|_u ,..., \frac{\partial}{\partial y^{n}}|_u ,\,u\in{\pi}^{-1}(x),\, x\in {\bf U} \}$ 
the local vertical distribution $\mathcal{V}_U$. The subspaces $\mathcal{V}_u $ and 
$\mathcal{H}_u$ are such that the following splitting of ${\bf
T}_u {\bf N}$ holds:
\begin{displaymath}
{\bf T}_u {\bf N}=\mathcal{V}_u \oplus \mathcal{H}_u ,\, \forall
\,\, u\in {\bf N}.
\end{displaymath}
This decomposition is invariant by the action of ${\bf GL}(n,{\bf R})$ and it defines a 
non-linear connection (a connection in the sense of Ehresmann([3])) on the principal fiber 
bundle ${\bf N}({\bf M},{\bf GL}(n,{\bf R}))$.

The local basis of the dual vector space ${\bf T^{*}_u N},\,u\in
{\bf N}$ is
\begin{displaymath}
\{ dx^{1}|_u ,...,dx^{n}|_u, \frac{{\delta}y^{1}}{F}|_u,...,
\frac{{\delta}y^{n}}{F}|_u \} ,
\end{displaymath}
\begin{equation}
\frac{{\delta}y^{i}}{F}|_u =\frac{1}{F}(dy^{i}+N^{i}_{j}dx^{j})|_u
, \quad i,j=1,...,n.
\end{equation}
This basis is dual to the basis (2.2.1).

\section{The Chern connection and other connections}
The non-linear connection defined above provides the possibility to
define an the Chern connection.
Let us consider $x\in{\bf M}$, $u\in{\bf T}_x {\bf M}\setminus\{0\}$ and ${\xi}\in{\bf T}_x 
{\bf M}$. We define the canonical projections
\begin{displaymath}
\pi :{\bf N} \longrightarrow {\bf M},\quad
 u\longrightarrow x,
\end{displaymath}
\begin{displaymath}
\pi _1:{\bf N}\times {\bf TM} \longrightarrow {\bf N},\quad
 (u,\xi _x)\longrightarrow u,
\end{displaymath}
\begin{displaymath}
\pi _2 :{\bf N}\times {\bf TM} \longrightarrow {\bf TM},\quad
(u,\xi _x)\longrightarrow \xi _x ,
\end{displaymath}
with $u\in \pi ^{-1}(x), x\in {\bf M}.$ Then the vector bundle $\pi^* {\bf TM}$ is 
completely determined as the minimal subset of ${\bf N}\times {\bf TM}$ by the equivalence relation 
defined in the following way: for every $u\in {\bf N} $ and
$(u,\xi) \in \pi^{-1} _1 (u)$,
\begin{displaymath}
(u,\xi)\in {\bf \pi^* TM}\quad \textrm{iff} \quad \pi \circ\pi
_2(u,\xi)=\pi\circ \pi_1(u,\xi).
\end{displaymath}

We can define in a similar way the vector bundle ${\bf \pi ^* SM}$ over ${\bf SM}$, being $\pi:{\bf 
SM}\longrightarrow {\bf M}$ the canonical projection in case of
positive homogeneous Finsler structures or ${\bf
PTM}:\longrightarrow {\bf M}$ in  case of absolutely homogeneous
Finsler structures.
\begin{definicion}
Let $({\bf M},F)$ be a Finsler structure. The fundamental and the Cartan tensors are defined 
in the natural local coordinate system $(x,y, {\bf U})$ by the
equations:
\begin{enumerate}
\item Fundamental tensor:
\begin{equation}
g(x,y):=\frac{1}{2}\frac{{\partial}^2 F^2(x,y) }{{\partial}y^i
{\partial}y^j }\, dx^i \otimes dx^j .
\end{equation}
\item Cartan tensor:
\begin{equation}
A(x,y):=\frac{F}{2} \frac{\partial g_{ij}}{\partial y^{k}}\,
\frac{{\delta}y^i}{F} \otimes dx^j \otimes dx^k=A_{ijk}\,
\frac{{\delta}y^i}{F} \otimes dx^j \otimes dx^k .
\end{equation}
\end{enumerate}
\end{definicion}
The Chern connection $\nabla$ is determined through the structure equations for the 
connection 1-forms as follows ([1]),
\begin{teorema}
Let $({\bf M},F)$ be a Finsler structure. The vector bundle ${\pi}^{*}{\bf TM}$ admits a 
unique linear connection characterized by the connection 1-forms $\{ {\omega}^i _j,\,\, 
i,j=1,...,n \} $ such that the following structure equations hold:
\begin{enumerate}
\item ``Torsion free'' condition,
\begin{equation}
 d(dx^{i})-dx^{j}\wedge
w^{i}_{j}=0,\quad i,j=1,...,n.
\end{equation}

\item Almost g-compatibility condition,
\begin{equation}
dg_{ij}-g_{kj}w^{k}_{i}-g_{ik}w^{k}_{j}=2A_{ijk}\frac{{\delta}y^{k}}{F},\quad
i,j,k=1,...,n.
\end{equation}
\end{enumerate}
\end{teorema}
The torsion freeness condition is equivalent to the absence of terms
containing $dy^{i}$ in the connection 1-forms ${\omega}^i _j $ and
also implies the symmetry of the connection coefficients $\Gamma^i
_{jk}$([2]):
\begin{equation}
w^{i}_{j}={\Gamma}^{i}_{jk}\, dx^{k},\quad
{\Gamma}^{i}_{jk}={\Gamma}^{i}_{kj}.
\end{equation}
The torsion freeness condition and almost $g$-compatibility determines the expression of the 
connection coefficients of the Chern connection in terms of the Cartan and fundamental 
tensor components ([1]).

Following theorem $(2.3.2)$, there is a coordinate-free
characterization of the Chern connection. Let us denote by $V( \tilde{X})$
the vertical and by $H(\tilde{X})$ the horizontal components (defined by the
non-linear connection on ${\bf N}$) of an arbitrary tangent vector
$\tilde{X}\in{\bf T}_u {\bf N} $. Then the following corollaries are
immediate consequences of {\it theorem 2.3.2}:
\begin{corolario}
Let $({\bf M}, F)$ be a Finsler structure. The almost $g$-compatibility condition $(2.3.4)$ is 
equivalent to the conditions:
\begin{equation}
{\nabla}_{V(\tilde{X})}g=2A(\tilde{X},\cdot,\cdot),
\end{equation}
\begin{equation}
 {\nabla}_{H(\tilde{X})} g=0, \forall \tilde{X}\in {\bf TN}.
\end{equation}
\end{corolario}
{\bf Proof}: using local natural coordinates and reading from {\it
theorem 2.3.4}, we have that the covariant derivative of the
metric is
\begin{displaymath}
\nabla(g)=(dg_{ij}-g_{kj}w^{k}_{i}-g_{ik}w^{k}_{j})\pi^*e^i\otimes
\pi^*e^j\,=\, 2A_{ijk}\frac{{\delta}y^{k}}{F}\otimes \pi^*
e^i\otimes\pi^* e^j.
\end{displaymath}
By the definition of covariant derivative along a direction and
nothing that $\, 2A_{ijk}\frac{{\delta}y^{k}}{F}$ is vertical, one
gets,
\begin{displaymath}
\nabla_{\tilde{X}} (g):=\, 2A_{ijk}\frac{{\delta}y^{k}}{F}(\tilde{X})\,\pi^*
e^i\otimes \pi^* e^j,\,\,\forall \tilde{X} \in {\bf TN}.
\end{displaymath}
From this formula follows the result.\hfill $\Box$
\begin{corolario}
Let $({\bf M}, F)$ be a Finsler structure. The torsion-free condition $(2.3.3)$ is equivalent 
to the following conditions:
\begin{enumerate}
\item Null vertical covariant derivative of sections of ${\pi}^*
{\bf TM}$: let $\tilde{X} \in {\bf TN}$ and $Y\in {\bf TM}$, then
\begin{equation}
{\nabla}_{{V(\tilde{X})}} {\pi}^* Y=0.
\end{equation}
\item  Let us consider $X,Y\in {\bf TM}$
and the associated horizontal vector fields $\tilde{X}=X^i\frac{\delta}{\delta x^i} $ and 
$\tilde{Y}=Y^i \frac{\delta}{\delta x^i} $. Then the following
equality holds:
\begin{equation}
{\nabla}_{\tilde{X}} {\pi}^* Y-{\nabla}_{\tilde{Y}}{\pi}^*
X-{\pi}^* ([X,Y])=0.
\end{equation}
\end{enumerate}
\end{corolario}
{\bf Proof}: as before we consider the torsion condition in local
coordinates. Then the local frame $\{e_j\}$ of $\Gamma{\bf TM}$
commutes, $[e_i,e_j]=0$. Using the symmetry in the connection
coefficients and the definition of the torsion operator (2.11),
one obtains that
\begin{displaymath}
{\nabla}_{\tilde{e_i}} {\pi}^* e_j-{\nabla}_{\tilde{e_j}}{\pi}^*
e_i-{\pi}^* ([e_i,e_j])={\nabla}_{\tilde{e_i}} {\pi}^*
e_j-{\nabla}_{\tilde{e_j}}{\pi}^*e_i\,=(\Gamma^k_{ij}-\Gamma^k_{ji})\,\pi^*e_k=0.
\end{displaymath}
In order to proof the second condition, it is as follows,
\begin{displaymath}
{\nabla}_{{\frac{\partial}{\partial y^i}}} {\pi}^* e_j:=\pi^* e_k\,w^k_j(\frac{\partial}{\partial y^i})=\pi^* e_k
\Gamma^k_{dj}dx^d(\frac{\partial}{\partial y^i})=0.\hfill \Box
\end{displaymath}

The curvature endomorphisms associated with the connection $w$ are
determined by the Cartan's second structure equations,
\begin{equation}
{\Omega}^{i}_{j}:=dw^{i}_{j}-w^{k}_{j}\wedge w^{i}_{k},\,
i,j,k=1,...,n.
\end{equation}
In local coordinates, these curvature endomorphisms are decomposed
in the following way,
\begin{equation}
{\Omega}^{i}_{j}=\frac{1}{2} R^{i}_{jkl}dx^{k}\wedge dx^{l}+
P^{i}_{jkl}dx^{k}\wedge
\frac{{\delta}y^{l}}{F}+\frac{1}{2}Q^{i}_{jkl}\frac{{\delta}y^{k}}{F}\wedge
\frac{{\delta}y^{l}}{F}.
\end{equation}
The quantities $R^i _{jkl}$, $P^i _{jkl}$ and $Q^i _{jkl}$ are called the hh, hv, and 
vv-curvature tensor components of the Chern connection. In
particular, for the Chern connections it holds that the last
component is identically zero $Q=0$ ([1,{\it chapter 3}]) for arbitrary Finsler
structures. The other tensors have the following expressions:
\begin{equation}
R^{i}_{jkl}=\frac{\delta \Gamma^i_{jk}}{\delta x^l}-\frac{\delta \Gamma^i_{jl}}{\delta x^k}- \Gamma^i_{hl}\Gamma^h_{jk}+\Gamma^i_{hk}\Gamma^h_{jl},
\end{equation}
\begin{equation}
P^i _{jkl}= -F\frac{\partial \Gamma^i_{jk}}{\partial y^l}.
\end{equation}

Let us mention other examples of linear connections that are
relevant in Finsler geometry([1]):
\begin{enumerate}
\item
 is Cartan's connection, which is metric compatible, but has non-trivial
torsion. It is determined by the following connection forms:
\begin{displaymath}
(\,^c\omega  )^k _{i} =\omega ^k _{i}\, +\, A^k _{i j}\frac{\delta
y^j}{F},\, i,j,k=1,...,n.
\end{displaymath}
\item Berwald's connection, defined by the $1$-form connection
forms
\begin{displaymath}
(\,^b\omega  )^k _{i} =\omega ^k _{i}\, +\, A^k _{i j}dx^k,\,
i,j,k=1,...,n.
\end{displaymath}
It is torsion-free, although it is not metric compatible.
\end{enumerate}

\chapter{Introduction to Berwald Spaces}
We will follow in this {\it chapter} the corresponding {\it
chapters} $10$ and $11$ of reference [1]. The proofs of the
following statements can be found in this reference.

\section{Definition and general properties of Berwald Spaces}
\begin{definicion}
A Finsler structure $F$ is said to be of Berwald type if the
coefficients of the Chern connection $\Gamma^i_{jk}$, written in
natural coordinates, do not depend on $y$.
\end{definicion}
There is a nice characterization of Berwald spaces: $({\bf M},F)$ is a Berwald space iff the Chern's
 connection leaves invariant the value of the finsler norm along any curve on {\bf M} (see for instance [1], [6] or [9]).

Berwald spaces are slightly different than Riemannian spaces, which
are contained in the finsler category. This make them more treatable
that other kind of Finsler spaces. Indeed, there is a complete
classification of Berwald spaces ([6]). From a physical point of
view, Berwald spaces can hold the Equivalence Principle, which lies
on the foundations of General Relativity.
\paragraph{}
A direct consequence of this definition is that for a Berwald
structure, the Chern connection defines {\it per se} a linear connection
on the manifold {\bf M}. Therefore, there is defined a covariant
derivative on {\bf M}:
\begin{displaymath}
\nabla_X W =\Big(\frac{d
W^i}{dt}|_{\sigma(t)}\,+W^j\,\Gamma^i_{jk}(\sigma(t))\Big)\frac{\partial}{\partial
x^i}|_{\sigma(t)},\quad T=\frac{d\sigma(t)}{dt}.
\end{displaymath}
There is the following result ([1]):
\begin{proposicion}
Let $({\bf M},F)$ be a Berwald space. Then:
\begin{enumerate}
\item Given any parallel vector field $W$ along a curve $\sigma$
in {\bf M}, its Finslerian norm $F(W)=\sqrt{g_{W}(W,W)}$  is
necessarily constant along $\sigma$.

\item For {\bf M} connected, its Minkowski linear spaces $({\bf
T}_x{\bf M},F_x)$ are all {linearly} isometric to each other.
\end{enumerate}
\end{proposicion}
There are several characterizations of Berwald spaces:
\begin{proposicion}
Let $(M,F)$ be a Finsler manifold. Then the following criteria are
equivalent:
\begin{enumerate}
\item The hv-curvature is vanishes: $P^i_{jkl}=0$.

\item The Cartan tensor is covariantly constant along all
horizontal directions on the slit tangent bundle ${\bf
TM}\setminus{0}$, $A_{ijk|l}$.

\item $({\bf M}, F)$ is a Berwald space.

\item $\big(\Gamma^i_{jk}y^jy^k\big)_{y^py^q}$ does not depend on
$y$.

\item $\big(\gamma^i_{jk}y^jy^k\big)_{y^py^q}$ does not depend on
$y$.
\end{enumerate}
\end{proposicion}

The following {\it proposition} also holds,
\begin{proposicion}
Let $({\bf M},F)$ be a Finsler structure. Then
\begin{enumerate}

\item The structure is Berwald.

\item $P^i_{jkl}=\,^bP^i_{jkl}=0,$ where $^bP^i_{jkl}$ is the
$hv$-curvature of the Berwald connection.

\item The $hh$-curvature is given by:
\begin{displaymath}
R^i_{jkl}=\frac{\partial \Gamma^i_{jl}}{\partial
x^k}-\frac{\partial \Gamma^i_{kl}}{\partial
x^j}-\Gamma^i_{hk}\Gamma^h_{jl}-\Gamma^i_{hl}\Gamma^h_{jk}.
\end{displaymath}
\end{enumerate}
\end{proposicion}
In order to write the following proposition, we need to have the
following definition
\begin{definicion}
A Finsler structure $({\bf M},F)$ is called locally Minkowski if
at each point $x\in {\bf M}$ there is a local coordinate system
such that $F(y)$ does not depend on $x$.
\end{definicion}
The locally Minkowski spaces are characterized by
\begin{proposicion}
Let $({\bf M},F)$ be a Finsler manifold. Then the following
statements are equivalent,
\begin{enumerate}
\item Both the $R$ and $P$ curvatures of the Chern connection
vanishes.

\item The structure is locally Minkowski.

\end{enumerate}

\end{proposicion}
A Finsler surface is a $2$-dimensional surface endowed with a
Finsler structure ([1, {\it chapter 4}]). In the case of surfaces,
the $hh$-curvature corresponds to the curvature scalar function
$K$, the analogous to the Gaussian curvature for Finsler geometry.
Similarly, the $hv$-curvature is defined by the Cartan invariant
$I$. However, for Berwald surfaces $I=0$.

In order to introduce the following result, note that the Flag
curvature at the point $x$ with flag $y$ and transverse edge $V$
is given by
\begin{displaymath}
K(x,y)=\frac{V^iy^jR_{ijkl}(x,y)y^lV^k}{g_{(x,y)}(V,V)g_{(x,y)}(y,y)-g^2_{(x,y)}(y,V)},
\end{displaymath}
where the evaluation of all the quantities is done at the point
$(x,y)$. One can now state Szab\'o rigidity theorem:
\begin{teorema}
Let $({\bf M},F)$ be a connected Berwald surface for the Finsler
function $F$ such that is strongly convex in all ${\bf
TM}\setminus\{0\}$. Then,
\begin{enumerate}
\item If the curvature $K=0$, then $F$ is locally Minkowski
everywhere.

\item If the flag curvature $K$ is not identically zero, then $F$
is Riemannian everywhere.
\end{enumerate}
\end{teorema}
This result restricts the existence of pure Berwald spaces (that
means, the ones which are not Riemannian or locally Minkowski) to
higher dimension than two.

\section{Examples of Berwald Spaces}
Following the end of the {\it section 3.1}, we give some examples
of Berwald spaces.
\begin{enumerate}

\item {\bf Riemannian Spaces in dimension $n$.}  They are characterized by the fact that the
fundamental tensor $g_{ij}$ defines a quadratic form on ${\bf M}$ given by the fundamental tensor. In particular, this implies that
This implies that $P^i_{jkl}=0$, which means that the space is
Berwald. Alternatively, one can check that in any natural
coordinate system, the connection coefficients of the Chern's
connection does not depend on $y$. Indeed, the connection
coefficients of the Chern connection of the Riemannian metric $g$
are equal to
the Christoffel symbols of the Levi-Civita connection.

\item {\bf Locally Minkowski Spaces in dimension $n$.} There is a natural coordinate
system where the fundamental tensor is constant on $x$. Therefore
the Christoffel ``type'' symbol $\gamma^i_{jk}=0$ as well as the
non-linear connection coefficients $N^i_j$ are zero (because they
linear combination of the $\gamma$ functions), in the given
natural coordinate system.

\item By Szab\'o's rigidity theorem ([6]), if we look for a $y$-global
Berwald structures which are not Riemannian or locally Minkowski, one needs to look for in dimensions
 higher dimensions than $2$. However, there are
Berwald local surfaces, as the following example due to Berwald
and Rund shows ([1, {\it section} 10.3]):

\begin{enumerate}

\item {\bf Example of a y-local Berwald Surface.}

In this example, {\bf M} is ${\bf R^2}.$ The Finsler function is
given by a function $\xi(x^1,x^2)$ that is a non-constant solution
of the PDE
\begin{displaymath}
\xi\frac{\partial \xi}{x^1}-\xi\frac{\partial \xi}{x^2}=0.
\end{displaymath}
The solutions are given implicitly by ([1] and references there)
\begin{displaymath}
x^1+x^2 \xi=\psi(\xi)
\end{displaymath}
where $\psi$ is an arbitrary analytic function of $xi$ such that
$\psi''\neq 0$. Finally, the Finsler function is
\begin{equation}
F(x,y)=y^2 (\xi +\frac{y^1}{y^2})
\end{equation}
$F<0$ if $y^2<0$. Therefore it is a $y$-local strong convex.

The Cartan invariant is $I=0$, while the sectional curvature is
\begin{displaymath}
K(x,y)=\frac{\psi''(\xi)}{(\xi+\frac{y^1}{y^2})^3(\psi'(\xi)-x^2)^3}.
\end{displaymath}
Therefore $K\neq 0$. From the form of the function $F$ one notes
that the structure is not analytical in the whole ${\bf T}_x{\bf
R}^2$.

\item {\bf Example of a $y$-global non-trivial Berwald space}.

In order to give an example of a $y$-global Berwald structure, we
use a Randers metric.

Our example is based on the
following result ([1, {\it section 11.6}]),
\begin{teorema}
Let $({\bf M},F)$ be a Randers Space. Denote the underlying
Riemannian metric by $a$, its Levi-Civita connection by
$\gamma^i_{jk}$ and the underlying $1$-form by $b$. Assume
\begin{enumerate}
\item $\parallel b\parallel_a \,<1.$

\item The covariant derivative respect the Levi-Civita connection
of b vanishes in all directions,
\begin{displaymath}
b_{j|k}:=\frac{\partial b_j}{x^k}-b_s\gamma^s_{jk}=0.
\end{displaymath}
Then the Randers space is of Berwald type. Conversely, if the
Randers space is of Berwald type, then above both conditions hold.
\end{enumerate}
\end{teorema}
{\bf Remark}. There is at least one topological restriction to the
above construction. The parallel condition is equivalent to the
existence of a global non-zero everywhere vector field. Therefore,
using Poincare-Hopf index theorem, for compact manifolds without boundary
surface the Euler characteristic $\chi{\bf M}$ must vanish.

The example that we present is the following. The base manifold is
given by ${\bf M}={\bf S}^2\times {\bf S}^1$. The Riemannian
metric is
\begin{displaymath}
a=(sin^2(\phi)d\theta\otimes d\theta\,+ d\phi\otimes
d\phi)\,+dt\otimes dt.
\end{displaymath}
The parallel $1$-form is given by
\begin{displaymath}
b(x,y)=\epsilon dt,\quad |\epsilon|<1.
\end{displaymath} In local coordinates a tangent vector $y\in {\bf T}_x{\bf M}$ can be written as
\begin{displaymath}
y=y^{\theta}\partial_{\theta} + y^{\Phi}\partial_{\Phi}+y^t\partial_t.
\end{displaymath}

Then, the Finsler function is
\begin{equation}
F(x,y)=\sqrt{sin^2(\phi)(y^{\theta})^2\,+(y^{\phi})^2\,+(y^t)^2}\,+\epsilon
y^t.
\end{equation}
Therefore, by {\it theorem} (3.2.1) this function $F$ defines a
Berwald structure on ${\bf S}^2\times {\bf S}^1$. It is clear that
this construction can be generalized to higher dimensions, with
similar constructions on ${\bf S}^n\times {\bf S}^1$.
\end{enumerate}
\end{enumerate}

\chapter{Review of the Theory of the Averaged Structures Associated with Finsler Structures}
This {\it chapter} follows quite closely {\it section 4} of [8],
where the original theory of averages of geometric structures was presented. As such, this
{\it chapter} does not constitute a original result of the present
memory, although it is fundamental for it. However, some of the
statements are proved in another way, while some of the proves
have been omitted for brevity of this memory.
\section{The Averaged of Linear Connections}
\begin{definicion}
Let $({\bf M},F)$ be a Finsler function. Then the indicatrix at
the point $x\in {\bf M}$ is the convex sub-manifold ${\bf I}_x\subset {\bf
t}_x{\bf M}$ defined by the condition that $(x,y)\in {\bf I}_x$
iff $F(x,y)=1$.
\end{definicion}
This is equivalent to the definition of
the tangent sphere in Riemannian Geometry. From this perspective,
a Finsler structure is a smooth collection $\{{\bf I}_x \subset {\bf
t}_x{\bf M},\, x\in {\bf M}\}$ of smooth,  convex tangent sets,
one at each point $x\in {\bf M}$, while a Riemannian structure is
a smooth collection of tangent ellipsoids.

Let $\tilde{X}$ be a tangent vector field along the horizontal
path
$\tilde{\gamma}:[0,1]\longrightarrow  {\bf N}$ 
connecting the points $u\in {\bf I}_x$ and $v\in {\bf I}_z$. The parallel transport 
associated with the Chern connection along ${\tilde{\gamma}}$ of a section  $S\in \pi^* {\bf 
TM}$ is denoted by $\tau _{{\tilde{\gamma}}}S$. The parallel transport along 
${\tilde{\gamma}}$ of the point $u\in {\bf I}_x$ is $\tau _{{\tilde{\gamma}}} 
(u)={\tilde{\gamma}} (1)\in \pi^{-1}(z)$. We say that $\tilde{\gamma}$ is horizontal
if the tangent vectors along $\tilde{\gamma}$ are horizontal. The horizontal lift of a path is defined using the 
non-linear connection defined on the bundle ${\bf
TN}\longrightarrow {\bf N}$.
\begin{proposicion}
(Invariance of the indicatrix by horizontal parallel transport)
Let $({\bf M}, F)$ be a Finsler structure, ${\tilde{\gamma}}:[0,1]\longrightarrow {\bf N}$ 
the horizontal lift of the path ${\gamma}:[0,1]\longrightarrow {\bf M}$ joining $x$ and $z$. 
Then the value of the function $F(x,y)$ is invariant along $\tilde{\gamma}$.
 In particular, let us consider ${\bf 
I}_x$ (resp (${\bf I}_z$)) to be the indicatrix over $x$ (resp
$z$). Then ${\tau}_{\tilde{\gamma}}({\bf I}_x )={\bf I}_z. $
\end{proposicion}
{\bf Proof:} Let $\tilde{X}$ be the horizontal lift in ${\bf TN}$
of the tangent
vector field $X$ along the path ${\gamma}\subset {\bf M}$ joining $x$ and $z$, $S_1, S_2 \in 
\pi^* ({\bf T}_x {\bf M})$. Then {\it corollary} $2.3.3$ implies ${\nabla}_{\tilde{X}} g (S_1 
,S_2)=2 A(\tilde{X},S_1 ,S_2)=0$ because the vector field $\tilde{X}$ is horizontal and the 
Cartan tensor is evaluated in the first argument. Therefore the value of the Finslerian norm 
$F(x,y)=\sqrt{g_{ij}(x,y)\, y^i \, y^j}$, $y\in {\bf T}_x {\bf M},\, Y$ with $Y=\pi^* y$ is 
conserved by horizontal parallel transport,
\begin{displaymath}
 {\nabla}_{\tilde{X}} (F^2(x,y))={\nabla}_{\tilde{X}} (g(x,y))(Y,Y)+2 
g(x,{\nabla}_{\tilde{X}} Y)=0,
\end{displaymath}
being $\tilde{X}\in {\bf TN}$ an horizontal vector. The first term is zero because the above 
calculation.
The second term is zero because of the definition of parallel transport of sections $\nabla 
_{\tilde{X}} Y=0$.
In particular the indicatrix ${\bf I}_x$ is mapped to ${\bf I}_z$ because parallel transport 
is a diffeomorphism. \hfill$\Box$

{\bf Remark}. Note the difference between this statement and the
statement of {\it proposition 3.1.2}: while {\it proposition}
4.1.2 applies to a general Finsler structure, {\it proposition}
3.1.2 refers to Berwald structures, where the Chern's connection
defines an affine connection on {\bf M} directly. Then, the parallel
transport along curves on {\bf M} makes sense.

\paragraph{}
Let us denote by $\pi^* _v {\bf \Gamma M}$ the fiber over $v\in
{\bf N}$ and by ${\bf \Gamma}_x {\bf M}$ the space of tensors
restricted to $x\in {\bf M}$; for every tensor $S_x \in {\bf
\Gamma} _x {\bf M}$ and $\, v\in \pi^{-1}(z)$, $z\in {\bf
U}\subset {\bf M}$ we consider the homomorphism:
\begin{displaymath}
\pi _2 | _v :{{\pi^*  _v {\bf \Gamma} {\bf M}}}\longrightarrow {\bf \Gamma} _z {{\bf  
M}},\quad S_v\longrightarrow S _z
\end{displaymath}
\begin{displaymath}
\pi ^*  _v :{\bf \Gamma} _z{{\bf M}}\longrightarrow {{\pi^*  _v
{\bf {\Gamma}M}}},\quad S _z\longrightarrow  \pi^* _v S_z .
\end{displaymath}
Let $S\in \,\Gamma\pi^* {\bf TM}$ and $S(u):=S_u\in\,
\pi^{-1}(u).$ Then, it holds that
\begin{equation}
\pi^*_u\,\pi_2|_v S(v)=S(u),\,\, u,v\in \pi^{-1}(x).
\end{equation}
If $S(v)\in \,\Gamma_v\,\pi^*\Gamma M$, $S(u)\in
\,\Gamma_u\,\pi^*\Gamma M$, the fibers over $u$ and $v$
respectively of the bundle $\pi^*{\bf TM}\longrightarrow {\bf n}$,
we have that in a local frame $S(v)=\xi^J(x)\,\pi^*_v e_J|_x$ and
respectively $S(u)=\xi^J(x)\,\pi^*_u e_J|_x$, where we are using
multi-index notation. Therefore,
\begin{displaymath}
\pi_2|_v S(v)=\pi_2\,\big(\xi^J(x)\,\pi^*_v
e_J|_x\big)=\xi^J(x)\,e_J|_x.
\end{displaymath}
Then,
\begin{displaymath}
\pi^*_u\,\pi_2|_v
S(v)=\pi^*_u\,\xi^J(x)\,e_J|_x\,=\xi^J(x)\,\pi^*_u e_J|_x=S(u).
\end{displaymath}
Note that for arbitrary $u,v\in \pi^{-1} (x)$ and an arbitrary element $S_u \in \pi^* _u{\bf 
TM}$, in general is not true that
\begin{displaymath}
\pi ^*  _u \pi _2 \mid _v :\pi^* _v{\bf \Gamma}  {{\bf  M}}\longrightarrow \pi^*  _u {{{\bf 
\Gamma}  {\bf M}}} ,\quad S _v \longrightarrow S _u
\end{displaymath}
because $\pi ^*  _u \pi _2 | _v S_v =\pi ^*  _u S_x$ and it is not the same than $S_u \in 
\pi^{-1}_1 (u)$, the evaluation of the section $S\in \pi^* \Gamma
M$ at the point $u$.
\begin{center}
***
\end{center}
We are now ready to define the
averaging operation,
\begin{definicion}
Consider the family of automorphims $A_w:=\{A_w:\pi ^*_w {\bf TM} \longrightarrow \pi^*_w {\bf 
TM}\}$ with $w\in \pi ^{-1}(x)$  with $x\in {\bf 
M}$. The average of this family of operators is  
the operator $A_x : {\bf T }_x {\bf M} \longrightarrow  {\bf T}_x{\bf M}$ such that:
\begin{displaymath}
<A_w\,>\,:=<\pi_2 |_u A\,\pi^* _u>_u S_x =\frac{1}{vol({\bf I}_x)}\big(\int _{{\bf I}_x} \pi _2 |_u  
A_u  \pi^* _u\, dvol \big)S_x ,
\end{displaymath}
\begin{equation}
vol({\bf I}_x)\:=\int _{{\bf I}_x} dvol ,\quad u\in
{\pi^{-1}(x)},\, S_x \in {\bf \Gamma }_x {\bf M};
\end{equation}
 $dvol$ is the standard volume form induced on
the indicatrix ${\bf I}_x$ from the Riemannian volume of the
Riemannian structure $({\bf T}_x{\bf M}\setminus\{0\},g_x )$.
\end{definicion}

{\bf Meaning of the Averaging Operation}

This definition of the averaged operation is new, compared with
other averages:
\begin{enumerate}
\item
 For instance, in the theory of characteristic
classes, integration along the fiber commutes with the exterior
differential and this is an essential point to prove Thom's
isomorphism theorem ([10]). The integration is in this example of
forms on fibers that are finite vectors spaces.

\item In Classical Mechanics, integration along the fiber is used
to derive a simplified averaged model, which in some circumstances
is simpler to analyze ([11]). This is also an integration along
fiber, where the fibers are invariant tori.

\end{enumerate}
In these both cases, the fiber bundle structure is similar: we
have a bundle $\pi:{\bf P}\longrightarrow {\bf M}$  and then we
calculate the integrals on $\pi^{-1}(x)$ for a given $x\in {\bf
M}$.

However, the averaging procedure that we propose is a bit more
involved. In our case, we have a double fiber structure:
\begin{displaymath}
\begin{array}{cccccc}
  & {\pi_1}  &  &  {\pi} & \\
\pi^*{\bf M} & \longrightarrow  &  {\bf N} &\longrightarrow &{\bf
M}.
\end{array}
\end{displaymath}
Although the composition is also a fiber bundle $\pi_1 \circ
\pi:\pi^*{\bf TM}\longrightarrow {\bf M}$, the integration that we
performed is on a lift of the fiber in the intermediate base
manifold ${\bf I}_x\subset {\bf N}$ on $\pi^*{\bf TM}$.  Therefore
we need in this case more structure that an ordinary fiber bundle
structure. In particular, we need to fix the lift.

Let $f\in {\bf \mathcal{F}M}$ be a real, smooth function on the manifold ${\bf M}$; ${\pi}^* 
f \in {\pi}^* {\bf \mathcal{F}}{\bf M}$ is defined in the
following way,
\begin{definicion}
Let $({\bf M}, F)$ be a Finsler structure, $\pi (u)=x$ and consider $f\in {\bf 
\mathcal{F}M}$. Then
\begin{equation}
\pi^* _u f=f(x),\quad \forall u\in\pi^{-1}(x).
\end{equation}
\end{definicion}
The definition is consistent because the function $\pi^* _v f$ is constant for every $v\in 
\pi^{-1}(x)$: $ \pi^* _u f=f(x)=\pi^* _v f,\, \forall u,v\in \pi^{-1}(x)$. Therefore the 
image is the constant value $f(x)$ for every $w\in \pi^{-1}(x)$.
${\pi}^* _u :{\bf T}_x {\bf M}\longrightarrow \pi ^* _u {\bf T M}$ is an isomorphism between 
${\bf T} _x {\bf M}$ and $\pi _1 ^{-1}(u)$, $\, \forall x\in {\bf
M}, \, u\in \pi^{-1}(x)$.

\paragraph{} Let us denote the horizontal lifting operator in the following way:
\begin{equation}
\iota:{\bf T}{\bf M}\longrightarrow {\bf T}{\bf N}, \quad
  X=X^i\frac{{\partial}}{{\partial}x^i}|_x \longrightarrow
\tilde{X}=X^i\frac{{\delta}}{{\delta}x^i}|_{u}=\iota(X),\, u\in
\pi^{-1}(x).
\end{equation}
This homomorphism is injective and the final result is a section of the tensor bundle ${\bf 
TN}$. In addition, it defines unambiguously a horizontal tangent vector $\tilde{X}\in 
\mathcal{H}_u$ for every tangent vector $X\in {\bf T}_x{\bf M}$. In our calculations,
We will also consider the 
restrictions of this map $\{\iota _u,\, u\in \pi^{-1}(x)\}$, such that 
$\iota _u (X)=(\iota X)_u $, $X\in {\bf T}_x{\bf M}$.
\paragraph{}

The following {\it proposition} is  the basis of the theory of the
averaged structures associated with Finsler structures. The
original proof can be found in reference [8], as well as all the
proofs of the results presented in this section,
\begin{teorema}
Let $({\bf M}, F)$ be a Finsler structure and $u\in \pi^{-1}(x)$, with $x\in {\bf M}$ and 
let us consider the respective Chern connection $\nabla$. Then for
each tangent field $X\in {\bf T}_x {\bf M}$ there is defined on
${\bf M}$ a covariant derivative  $\tilde{{\nabla}}_X $ such that
\begin{enumerate}
\item $\forall  X\in {\bf T}_x {\bf M}$ and $Y\in  {\bf TM}$ the covariant derivative of $Y$ 
in the direction $X$ is given by the following average:
\begin{equation}
\tilde{{\nabla}}_X Y =<\pi _2|_u {\nabla}_{{\iota}_u (X)} {\pi}^* _v Y\,>_u 
,\,\, u\in {\bf I}_x\subset \pi^{-1}(x)\subset {\bf N} .
\end{equation}

\item For every smooth function $f\in {\bf \mathcal{F}}{\bf M}$ the covariant derivative is 
given by the following average:
\begin{equation}
\tilde{{\nabla}}_X f =<\pi _2 |_u \nabla _{\iota _u (X)} \pi^* _v
f>_u =X\cdot(f).
\end{equation}
\end{enumerate}
\end{teorema}
{\bf Proof:} The argument follows in the following way. Consider
the convex sum of linear connections $t_1\nabla_1+...
t_p\nabla_p$, $t_1+...+t_p=1$; the connections are linear
connections on {\bf M}. It is well known that $t_1\nabla_1+...
t_p\nabla_p$ is also a linear connection. Now, consider the
manifold ${\bf \Sigma}_x\subset \pi^{-1}(x)\subset {\bf N}$ and a
set of connections on {\bf M}, all of them labelled by points on
${\Sigma}$, so there is a map $\Theta:{\bf M}\longrightarrow ({\bf
R}^+)$ such that $\int_{{\bf \Sigma}_x} \Theta=1$ and that
$\Theta\geq 0$. Then, one can use a limit argument (work in
progress) to show that the averaged of the family of connections
$\{\nabla_u\}$ defines also a linear connection on {\bf M}. To
apply to our case this argument, we only need to specify that
$\Sigma_x={\bf I}_x$ and that $\Theta(u)=dvol\,\pi_2|_u\nabla_{\iota_u}
\pi^*$, where the right hand side must be understood for fixed
$u\in {\bf I}_x$ and as acting on sections of $\Gamma {\bf
M}$.\hfill$\Box$

\paragraph{}
The averaged covariant derivative commutes with contractions:
\begin{displaymath}
\tilde{\nabla}_X [\alpha (Z)]=<\pi_2|_u\nabla_{\iota_u (X)}\pi^*
(\alpha (Z))>_u:=\iota_{\tilde{\nabla}_X (Z)}\alpha
\,+\iota_{Z}\tilde{\nabla}_X \alpha.
\end{displaymath}
The extension of the covariant derivative $\tilde{\nabla}_X$ acting on sections of ${\bf 
\Gamma}^{(p,q)}{\bf M}$ is performed in the usual way,
\begin{displaymath}
\tilde{\nabla}_X
\,K(X_1,...,X_s,\alpha^1,...,\alpha^r)=\tilde{\nabla}_X\,
K(X_1,...,X_s,\alpha^1,...,\alpha^r)-
\end{displaymath}
\begin{displaymath}
-\sum^s_{i=1}K(X_1,...,\tilde{\nabla}_X
X_i,...,X_s,\alpha^1,...,\alpha^s) +\sum^r_{j=1}
K(X_1,...,X_s,\alpha^1,...,\tilde{\nabla}_X \alpha^j,...\alpha^r).
\end{displaymath}

We denote the affine connection associated with the above covariant derivative by 
$\tilde{{\nabla}}$: for every section $Y\in {\bf TM}$, $\tilde{\nabla}Y\in {\bf T}^*_x{\bf  
M}\otimes  {\bf TM},\, x\in {\bf M}$ is given by the action on pairs $(X,Y)\in {\bf T 
}_x{\bf  M}\otimes  {\bf TM}$,
\begin{equation}
\tilde{\nabla}(X,Y):=\tilde{\nabla}_X Y.
\end{equation}

{\bf Remark}. From the proof of {\it theorem} 4.1.5 one easily
recognize that the result can be applied to any other linear
connection defined in the bundle $\pi^* {\bf TM}$.

Let us calculate the torsion of the connection $\tilde{\nabla}$.
Then the torsion is given for arbitrary vector fields $X,Y \in
{\bf  TM}$ by
\begin{displaymath}
T_{\tilde{\nabla}}(X,Y)=<\pi _2 |_u \nabla _{\iota _u (X)} \pi^* _w >_u Y -<\pi 
_2 |_u  \nabla _{\iota _u (Y)} \pi^* _w >_u X -[X,Y]=
\end{displaymath}
\begin{displaymath}
=<\pi _2 |_u \nabla _{\iota _u (X)}  \pi^* _u >_u Y -<\pi _2 |_u 
\nabla _{\iota _u (Y)} \pi^* _u >_u X
\end{displaymath}
\begin{displaymath}
 -<\pi _2 |_u  \pi^* _u [X,Y]>=
\end{displaymath}
\begin{displaymath}
=<\pi _2 |_u  \big(\nabla _{\iota _u (X)}\pi^* Y -\nabla _{\iota _u (Y)}\pi^*X 
-\pi^*[X,Y] \big) >_u \,=0,
\end{displaymath}
because the torsion-free condition of the Chern connection.
Therefore,
\begin{proposicion}
Let ({\bf M},F) be a Finsler structure with averaged connection $\tilde{\nabla}$. Then the 
torsion $T_{\tilde{\nabla}}$ of the average connection obtained from the Chern connection is 
zero.
\end{proposicion}
\begin{center}
***
\end{center}
Let us consider the following (non-degenerate) tensors,
\begin{displaymath}
 g_t =(1-t)g +t h, \,t\in[0,1].
\end{displaymath}
$g_t $ defines a Finsler structure in ${\bf M}$. The associated Chern's connection are 
denoted by $\nabla _t $. In a similar as in {\it theorem 4.1.5}
the following result is proved:
\begin{teorema}
Let $({\bf M},F)$ be a Finsler manifold and $ g_t =(1-t)g +th $,
$t\in[0,1]$. Then the operator
\begin{equation}
\tilde{\nabla} _t = \frac{1}{vol({\bf I}_x)}\int_{{\bf I}_x} \pi _2 |_u  
\nabla  \pi^* _v
\end{equation}
is a linear connection on ${\bf M}$ with zero torsion for every
$t\in[0,1]$.
\end{teorema}

\section{Structural theorems}

We consider some results obtained in [8] relating
geometric objects of the averaged connection and their related
averaged objects. The results can be applied to the averaged
connection of any linear connection on $\pi^*{\bf TM}$.

If $\mathcal{E}\longrightarrow {\bf N}$ is an arbitrary vector bundle over {\bf N}, for a 
given parallel transport $\tau $ along an arbitrary path $x_t$ with tangent vector 
$\dot {x}_t\in {\bf T}_u{\bf N}$, the covariant derivative of a
section $S$ is given by the expression
\begin{displaymath}
\nabla _{\dot{x}_t} S =\lim _{\delta \rightarrow 0}\frac{1}{{\delta}}\big( \tau ^{t+\delta} 
_t S (x_{t+\delta}) -S (x_t)\big).
\end{displaymath}
Applying this formula to the Chern connection,
\begin{displaymath}
\tilde{\nabla}_{X} S= <\pi _2 |_u(t)  \lim _{\delta\rightarrow 0} \frac{1}{\delta}\big( \tau 
^{t+\delta} _t\pi^* _{u(t+\delta)} S(x_{t+\delta}) -\pi^* _{u(t)}
S(x_t )\big)>_{u(t)},
\end{displaymath}
with $u(t+\delta)\in {\pi^{-1}(x(t+\delta))}$. Interchanging the limit and the average 
operation (this can be done, because both integrals are performed
on the same manifold) one obtains,
\begin{displaymath}
\tilde{\nabla}_X S =\lim _{\delta \rightarrow 0} \frac{1}{\delta}<\pi _2 |_u \big( \tau 
^{t+\delta} _t \pi^* _{u(t+\delta)} S(x_{t+\delta}) -\pi^*
_{u(t)}S(x_t )\big)>_{u(t)}.
\end{displaymath}
This interchange can be done because the integration is performed in ${\bf I}_x$, not 
depending of the limit label $\delta$. Therefore,
\begin{displaymath}
<\pi _2 |_u \pi^* _{u(t)} S(x_t )>_{u(t)} =<\pi _2 |_u \pi^* _{u(t)} S^{\mu}(x_t 
)\frac{\partial}{\partial x^{\mu}}>_{u(t)}=
\end{displaymath}
\begin{displaymath}
=S^{\mu}(x_t ) <\pi _2 |_u \pi^* _{u(t)} \frac{\partial}{\partial x^{\mu}}>_{u(t)}=S^{\mu} 
\frac{\partial}{\partial x^{\mu}}.
\end{displaymath}
Then one can conclude that the expression
\begin{displaymath}
<\pi _2 |_u \tau ^{t+\delta} _t \pi^* _{u(t+\delta)} >_u
\end{displaymath}
plays the role of the parallel transport operation for the average connection 
$\tilde{\nabla}$,
\begin{teorema}
Let $({\bf M},F)$ be a Finsler structure with associated Chern's connection $\nabla $ and 
with average connection $\tilde{\nabla}$. Then the parallel transport associated with 
$\tilde{\nabla}$ along a short path of parameter length $\delta t$ is given by
\begin{equation}
(\tilde{\tau}^{t+\delta}_t )_{x_t}S\,:=<\pi _2 |_u \tau ^{t+\delta} _t \pi^* _{u(t+\delta)} 
S(x_{t+\delta})>_{u(t)} ,\quad  S_x \in {\bf \Gamma} _{x_{t}}{\bf
M}.
\end{equation}
\end{teorema}
 {\bf Proof:} It is immediate from the definition of the covariant derivative in terms of 
infinitesimal parallel transport; let us define the section along
$\gamma _t$ by
\begin{displaymath}
\tilde{\tau}S(x_{t+\delta})=\tilde{\tau}^t _{t+\delta}S(x_t ),
\end{displaymath}
that is the parallel transported value of the section $S$ from the point $x_{t+\delta}$ to 
the point $x_{t}$. Then it follows from the general definition of
covariant derivative that
\begin{displaymath}
\tilde{\nabla}_X (\tilde{\tau} S )=\lim _{\delta \rightarrow 0}\frac{1}{\delta} 
\big(\tilde{\tau}^{t+\delta} _t S(x_{t+\delta} )-S(x_t)\big)=0.
\end{displaymath}
When this condition is written in local coordinates, it is equivalent to a system of ODEs 
and the result is obtained from uniqueness and existence of
solutions of ODEs.\hfill$\Box$

However, in order to check the consistency of this definition, we should check that the 
composition rule holds:
\begin{displaymath}
\tilde{\tau}^{t+\delta} _t \circ \tilde{\tau}^{t+2\delta} _{t+\delta} 
=\tilde{\tau}^{t+2\delta} _t .
\end{displaymath}
The proof is directly obtained calculating the above composition. This calculation reveals 
the reason to take $\tilde{\nabla}_1$ as the average connection if one wants to preserve the 
rule of ``average parallel transport as the parallel transport of the average connection" 
for finite paths:
\begin{displaymath}
\tilde{\tau}^{t+\delta} _{t}\circ \tilde{\tau}^{t+2\delta} _{t+\delta} =\frac{1}{vol({\bf 
I}_{u(t)})}\frac{1}{vol({\bf I}_{u(t+\delta)})}\int _{{\bf I}_{x_t}} \int _{{\bf 
I}_{x_{t+\delta}}} \pi _2 |_{u(t)} {\tau}^{t+\delta} _{t}\pi^*_{u(t+\delta)}\pi _2 
|_{u(t+\delta)} \circ \,
\end{displaymath}
\begin{displaymath}
{\tau}^{t+2\delta} _{t+\delta}\pi^* _{u(t+2\delta)}=\frac{1}{vol({\bf 
I}_{u(t)})}\frac{1}{vol({\bf I}_{u(t+\delta)})}\int _{{\bf I}_{x_t}} \int _{{\bf 
I}_{x_{t+\delta}}} \pi _2 |_{u(t)} {\tau}^{t+\delta} _{t}\circ {\tau}^{t+2\delta} 
_{t+\delta}\pi^* _{u(t+2\delta)}=
\end{displaymath}
\begin{displaymath}
\frac{1}{vol({\bf I}_{u(t)})}\int _{{\bf I}_{x_t}}\pi _2 |_{u(t)} {\tau}^{t+2\delta} 
_{t}\pi^* _{u(t+2\delta)}=\tilde{\tau}^{t+2\delta} _t .
\end{displaymath}
\paragraph{}
We use a well known formula in order to express curvature endomorphisms as an
infinitesimal parallel transport ([1]): denote by $\gamma _t
:[0,1]\longrightarrow {\bf M}$ the infinitesimal parallelogram
built up from the vectors $X,Y \in {\bf T}_x {\bf M}$ with lengths
equal to $\delta t$ constructed using parallel transport along the
integral curves of $X, Y,-X,-Y$ through a short time $\delta t$
and where length is measure using the Finslerian length. Then for
every linear connection, the curvature endomorphisms are given by
the formula
\begin{equation}
\Omega(X,Y)=-\frac{d\tau (\gamma _t )}{dt}|_{t=0}.
\end{equation}
It can be written formally like
\begin{equation}
I+\delta t \Omega(X,Y)=-\tau _{d\tilde{\gamma}}.
\end{equation}
Let us denote by $\tilde{\Omega}=\tilde{R}:=R^{\tilde{\nabla}}$ the curvature of 
$\tilde{\nabla}$ and let us recall the hh-curvature of the Chern
connection,
\begin{teorema}
Let $({\bf M}, F)$ be a Finsler structure. Let $\iota({X_1}), \iota({X_2})$ be the 
horizontal lifts in ${\bf T}_u {\bf N}$ of the linear independent vectors $X_1, X_2 \in {\bf 
T}_x {\bf M}$. Then for every section $Y\in {\bf \Gamma}{\bf M}$,
\begin{equation}
\tilde{R}_x (X_1 ,X_2 )Y= <\pi _2 R_u (\iota _u({X_1}), \iota _u ({X_2}))  
\pi^*_u  Y>_u ,\, u\in {\bf I}_x\subset \pi^{-1}\subset {\bf N}.
\end{equation}
\end{teorema}
{\bf Algebraic Proof.} Let us assume a local frame of vector
fields. Then we can write the value of the averaged curvature
endomorphism as
\begin{displaymath}
\tilde{R}_{u}(X,Y)Z=<\pi_2|_u\nabla_{\iota_u (X)}\pi^*|_u
\cdot<\pi_2|_v\nabla_{\iota_v (Y)}\pi^*|_v Z>>-
\end{displaymath}
\begin{displaymath}
 -<\pi_2|_u\nabla_{\iota_u (Y)}\pi^*|_u
\cdot<\pi_2|_v\nabla_{\iota_v (X)}\pi^*|_v Z>>-
\end{displaymath}
\begin{displaymath}
-<\pi_2|_u\nabla_{\iota_u ([X,Y])}\pi^*|_u Z>,\,\,X,Y,Z\in
\Gamma{\bf TM}.
\end{displaymath}
Using the relation (4.1.1) one can reduce the above double
integral to a single integral. For instance,
\begin{displaymath}
<\pi_2|_u\nabla_{\iota_u (X)}\pi^*|_u
\cdot<\pi_2|_v\nabla_{\iota_v (Y)}\pi^*|_v
Z>>=<\pi_2|_u\nabla_{\iota_u (X)}\nabla_{\iota_v (Y)}\pi^*|_u Z>
\end{displaymath}
Therefore,
\begin{displaymath}
\tilde{R}_{u}(X,Y)Z\,=<\pi_2|_u\nabla_{\iota_u (X)}\nabla_{\iota_v
(Y)}\pi^*|_u Z>\,-<\pi_2|_u\nabla_{\iota_u (Y)}\nabla_{\iota_v
(X)}\pi^*|_u Z>-
\end{displaymath}
\begin{displaymath}
-<\pi_2|_u\nabla_{\iota_u ([X,Y])}\pi^*|_u Z>=<\pi_2|_u R_u
([X,Y])\pi^*|_u Z>:=<R>_x(X,Y)Z.
\end{displaymath}\hfill $\Box$

From this second proof one can think that given two averaged
objects, if we multiply them, the product of averages is the
average of the product. However this is not true, as the following
counterexample shows,
\begin{displaymath}
\tilde{\nabla}<g>=<\pi_2|_u\nabla_{{\iota}_u X}\pi^*_u<\pi_2|_v
g_{ij}(x,v)\pi^*_ve^i\otimes \pi^*_ve^j>\,>\neq
\end{displaymath}
\begin{displaymath}
\neg<\pi_2|_u\nabla_{{\iota}_u X} g_{ij}(x,u)\pi^*_ve^i\otimes
\pi^*_u e^j>>
\end{displaymath}
because the coefficients $g_{ij}$ live on ${\bf I}_x$ and not on
{\bf M}.

\chapter{Some Applications of the Averaged Connection}
In this {\it chapter} we present some applications of the
averaged connection.
\section{Metric compatibility of the Averaged Connection}
It is interesting to consider when the averaged connection obtained from the Chern connection 
is Riemann metrizable, that means, when exists a Riemannian metric $\tilde{h}$ such that 
$\tilde{\nabla}=\nabla^{\tilde{h}}$. The basic result is the
following
\begin{proposicion}
Let $({\bf M},F)$ be a Finsler structure. Then the averaged connection ${\tilde{\nabla}}$ of 
the Chern connection ${\nabla}$ is a metric irreducible connection iff the Holonomy 
group $Hol({\tilde{\nabla}})$ is a Berger group.
\end{proposicion}
{\bf Proof:} Suppose that $\tilde{\nabla}$ is metrizable. Then there is a Riemanian metric such 
that $\tilde{\nabla}=\nabla^{\tilde{h}}$, that is, the Levi-Civita
connection of $h$. Since the torsion $T_{\tilde{\nabla}}=0$, it
implies $\tilde{\nabla}$ is a Riemannian connection and therefore
in the case of irreducible metrics, the holonomy group
$\tilde{\nabla}$ is a Berger group.

Conversely, let us suppose that $Hol(\tilde{\nabla})$ is an
irreducible Berger group. Then it is compact. Then we can define
the operation:
\begin{displaymath}
\int _{Hol(\tilde{\nabla})}d\tau;\,\,\,\quad \tau\in
Hol(\tilde{\nabla}).
\end{displaymath}
$d\tau$ is an invariant Haar measure of the Berger group $Hol(\tilde{\nabla})$. In particular we can use the 
Szabo's construction in [6] to define the following scalar product
on ${\bf T}_x {\bf M}$:
\begin{equation}
\tilde{h}_x(X,Y) =\int _{Hol(\tilde{\nabla})}(\tau^* X,\tau^* Y)^* \,d\tau;\, X,Y \in {\bf 
T}_x {\bf M}.
\end{equation}
$(,)^*$ is an arbitrary scalar product on ${\bf T}_x {\bf M}$. One extends this scalar product 
to the whole manifold using the holonomy group, defining a Riemannian metric $\tilde{h}$ 
that is conserved by $\tilde{\nabla}$.\hfill$\Box$

\section{Geodesic Equivalence Problem}

In order to clarify the relation between $h$ and $\tilde{h}$, we
use the notion of geodesic rigidity to obtain a partial answer to
this question.
\begin{definicion}
Two Riemannian metrics $h$ and $\bar{h}$ living on the manifold ${\bf M}$ with $dim({\bf 
M})\,\geq 2$ are geodesically equivalent if their sets of un-parameterized geodesics 
coincide.
The manifold {\bf M} is called geodesically rigid if every two geodesically equivalent 
metrics are proportional.
\end{definicion}

\begin{corolario}
Under the above hypothesis than before, $h$ and $\tilde{h}$ have
the same Levi-Civita connection. Therefore $h$ and $\tilde{h}$ are
geodesically equivalent.
\end{corolario}
{\bf Proof:} By definition $h$ is the Levi-Civita connection of
$\tilde{\nabla}$. On the other hand,
\begin{displaymath}
\big(\tilde{\nabla}_Z(\tilde{h})\big)(X,Y)=0,\quad \forall X,Y,Z
\in \Gamma {\bf TM}.
\end{displaymath}
because it has been extended using the holonomy group. Therefore
because $\tilde{\nabla}$ is also torsion three, it is the
Levi-Civita connection of $\tilde{h}$. \hfill$\Box$

{\bf Remark}. The above {\it corollary} is stronger than
geodesically equivalence condition between two metrics, because
the connection is already determine.

\begin{corolario}
Let $({\bf M}, F)$ be a Finsler structure such that ${\bf M}$ is Riemannian geodesically rigid, 
with $Hol(\tilde{\nabla})$ a Berger group. Then 
$h=C\tilde{h}$.
\end{corolario}

\paragraph{}
Matveev solved the problem of geodesically rigidity in Riemannian manifolds 
(see for instance ref. [12], [13] and [14]): to decide wether or
not two given metrics with the same geodesics are equivalent. In
particular, for hyperbolic manifolds, being Riemannian
geodesically rigid, one obtains
\begin{corolario}
Let $({\bf M}, F)$ be a Finsler structure such that ${\bf M}$ is a
closed manifold and
such that $\tilde{h}$ is an 
hyperbolic metric such that $Hol(\tilde{\nabla})$ is 
a Berger group. Then $h=C\tilde{h}$.
\end{corolario}

For a Berwald space the Holonomy group 
$Hol({\nabla})$ is compact. Then the holonomy group $Hol(\tilde{\nabla})$ 
is also compact and is a Berger group. Therefore it is a direct
consequence from a theorem of Matveev ([14]) the following
\begin{corolario}
Let $({\bf M},F)$ be a Berwald structure such that ${\bf M}$ admits an hyperbolic Riemannian metric. Then 
$h=\tilde{h}$ and $\tilde{\nabla}h =0,$ where $h$ and $\tilde{h}$
are defined as before.
\end{corolario}

\paragraph{}
In Finsler Geometry, the Finsler function is in general not reversible $(F(x,y)\neq F(x,-y))$.
Therefore it has sense the notion of geodesic reversibility.
Consider the following piece-wise differentiable curve $\gamma
\cup \beta$ where $\gamma(s),\,s\in [0,s_0]$ is a geodesic of the
Chern connection and $\beta$ is a geodesic but that start at the
end of $\gamma$ and has reverted the final vector of the first
geodesic. Let us close with another simple curve that is simple $\Delta$
the above curve. We call it closed almost-geodesic triangle.
\begin{definicion}
We say that the structure $({\bf M},F)$ is geodesically reversible
if every almost-geodesic triangle is retractible. Otherwise is
geodesic irreversible.
\end{definicion}
Examples of geodesically reversible structure are Riemannian
manifolds and reversible Finsler manifolds. A non-trivial example
is provided by Randers structures. One obtains the following {\it
proposition},
\begin{proposicion}
Let $({\bf M},F)$ be a Berwald structure. Then,
\begin{enumerate}
\item It is geodesically rigid if the Szabo's metric is Riemannian
geodesically equivalent. In this case, there is a Riemannian
metric with the same geodesics.

\item It is geodesically reversible.
\end{enumerate}
\end{proposicion}
A similar result is obtained by $Matveev$ ({\it theorem 1} of
[9]).

\section{A rigidity property for Berwald Spaces}

We start considering a generalization of some well known
properties of linear connections over ${\bf M}$ ([3], section 5.4)
to linear connections defined on the bundle $\pi ^* {\bf
TM}\rightarrow {\bf N}$.

Given two linear connections $^1\nabla$ and $^2\nabla$ on the
bundle $\pi ^* {\bf TM}\rightarrow {\bf N}$, the difference
operator
\begin{displaymath}
B:\Gamma {\bf TM}\otimes \Gamma {\bf TM}\rightarrow \pi^*{\bf
\Gamma TM}
\end{displaymath}
\begin{displaymath}
B_u(X,  Y)=\,^1\nabla_{\iota_u(X)}\pi^*_u Y
-\,^2\nabla_{\iota_u(X)}\pi^*Y,
\end{displaymath}
\begin{displaymath}
 u\in {\bf N}, \, X,Y\in {\bf \Gamma TM}
\end{displaymath}
is an homomorphism that holds the Leibnitz rule on $Y$ and it is $\mathcal{F}$-linear on $X$.

The symmetric and skew-symmetric components $S$ and $A$ of B are
defined in the following way
\begin{displaymath}
S_u\,:{\bf \Gamma TM}\times {\bf \Gamma TM}\longrightarrow \pi^*_u
{\bf TM}
\end{displaymath}
\begin{displaymath}
S_u(X,Y):=\frac{1}{2}\,\big(B_u (X,Y)+B_u( Y,X)\big).
\end{displaymath}
\begin{displaymath}
u\in \pi^{-1}(x),\quad X, Y\in {\bf \Gamma TM}.
\end{displaymath}
Then, the following relation holds for arbitrary vector fields
$X,Y\in \Gamma{\bf TM}$,
\begin{displaymath}
2S_u(X,Y)=\,^1\nabla{(\iota_u (X))}\pi^*_u Y -\,^2\nabla_{(\iota_u
(X))}\pi^*_u Y +(\,^1\nabla_{(\iota_u (Y))}\pi^*_u X
-\,^2\nabla_{(\iota_u (Y))}\pi^*_u X)=
\end{displaymath}
\begin{displaymath}
=\big(\,^1\nabla_{(\iota_u (X))}\pi^*_u Y\,+\,^1\nabla_{(\iota_u
(Y))}\pi^*_u X\big)-\,\big(\,^2\nabla_{(\iota_u (X))}\pi^*_u
Y\,+\,^2\nabla_{(\iota_u (Y))}\pi^*_u X\big).
\end{displaymath}
The skew-symmetric part $A$ is defined in a similar way,
\begin{displaymath}
A_u\,:{\bf \Gamma TM}\times {\bf \Gamma TM}\longrightarrow \pi^*_u
{\bf TM}
\end{displaymath}
\begin{displaymath}
A_u(X,Y):=\frac{1}{2}\,\big(B_u(X,Y)-B_u(Y, X)\big),
\end{displaymath}
\begin{displaymath}
\forall u\in \pi^{-1}(x),\quad X\in {\bf T}_x{\bf M}\quad ,Y\in
{\bf \Gamma TM}.
\end{displaymath}
As for the torsion, one can define the symmetric and
skew-symmetric parts $S$ and $A$ as a family of operators, because
the above definitions are point-wise.

Then, the following relation holds for arbitrary vector fields
$X,Y\in \Gamma{\bf TM}$,
\begin{displaymath}
2A_u(X,Y)=\nabla_{1(\iota_u (X))}\pi^*_u Y -\nabla_{2(\iota_u
(X))}\pi^*_u Y -(\nabla_{1(\iota_u (Y))}\pi^*_u X
-\nabla_{2(\iota_u (Y))}\pi^*_u X)=
\end{displaymath}
\begin{displaymath}
=Tor_u(\nabla_1)(X,Y) -Tor_u(\nabla_2)(X,Y).
\end{displaymath}
Since this relation holds point-wise for all $u\in \pi^{-1}(x)\in
\subset {\bf N}$ we can write
\begin{equation}
2A(X,Y)=Tor(\nabla_1)(X, Y)-Tor(\nabla_2)(X,Y).
\end{equation}
\begin{definicion}
Let $\nabla$ be a linear connection on the vector bundle
$\pi^*{\bf TM}\longrightarrow {\bf N}$ with connection
coefficients $\Gamma^i _{jk}$. The geodesics of $\nabla$ are the
parameterized curves $x:[a,b]\longrightarrow {\bf M}$ solutions of
the differential equations
\begin{equation}
\frac{d^2 x^i}{ds^2}+\Gamma^i _{jk}(x,\frac{dx}{ds})
\frac{dx^j}{ds}\frac{dx^k}{ds}=0,\quad i,j,k=1,...,n,
\end{equation}
where $\Gamma^i _{jk}(x,y)$ are the connection coefficients of
$\nabla$.
\end{definicion}
This differential equation can be written as
\begin{equation}
\nabla_{\iota_u (X)}\pi^*_u X=0,\quad u=\frac{dx}{ds}
\end{equation}
\paragraph{}

The following propositions are direct generalizations of the
analogous results for affine connections ([3]).
\begin{proposicion}
Let $^1\nabla$ and $^2\nabla$ be linear connections on the vector
bundle ${\pi^*}{\bf TM}\rightarrow {\bf N}$. Then the following
conditions are equivalent:
\begin{enumerate}
\item The connections $^1\nabla$ and $^2\nabla$ have the same
geodesic curves on {\bf M},

\item $B_u(X,X)=0$,

\item $S_u=0$,

\item $B_u=A_u,\,\forall u\in {\bf N}.$
\end{enumerate}
\end{proposicion}
{\bf Proof}. The proof follows the lines of ref. [3, pg 64-65]:
\begin{enumerate}
\item $({\bf a\Rightarrow b})$. If $^1\nabla$ and $^2\nabla$ have
the same geodesics, then they have the same geodesic equations.
Therefore
\begin{displaymath}
^1\nabla_{\iota_u(X)} \pi^*(X) =0\,\,\,\Leftrightarrow
\,\,\,^2\nabla_{\iota_u(X)} \pi^*(X) =0.
\end{displaymath}
This is implies that
\begin{displaymath}
B_u(X,X)=\,^1\nabla_{\iota_u(X)} \pi^*(X) -^2\nabla_{\iota_u(X)}
\pi^*(X)=0.
\end{displaymath}

 \item $({\bf b\Rightarrow c})$. It is consequence of linearity,
\begin{displaymath}
0=2B_u (X+Y,X+Y)=2S_u(X,Y).
\end{displaymath}

\item $({\bf c\Rightarrow d})$. Trivially from the definition of
$B$, $S$ and $A$.

\item $({\bf d\Rightarrow a})$. If $B=A$, implies $S=0$. In
particular $S_u(X,X)$=0, which implies
\begin{displaymath}
^1\nabla_{\iota_u(X)} \pi^*(X) = \,^2\nabla_{\iota_u(X)} \pi^*(X).
\end{displaymath}
From this relation and from existence and uniqueness of solutions,
the parameterized geodesics of $\,^1\nabla$ and $\,^\nabla$
coincide.\hfill $\Box$
\end{enumerate}
\begin{proposicion}
Let $^1\nabla$ and $^2\nabla$ be linear connections on the vector
bundle $\pi ^* {\bf TM}\rightarrow {\bf N}$ such that they have
the same covariant derivative along vertical directions. Then
$^1\nabla=\,^2\nabla$ iff they have the same parameterized
geodesics and $Tor(\,^1\nabla)=Tor(\,^2\nabla)$.
\end{proposicion}
{\bf Proof:} If $^1\nabla =\, ^2\nabla$, then they have the same
parameterized geodesics and torsion tensors. Conversely, if the
geodesics are the same, the torsion is the same, then $B=0$. Since
by hypothesis both connections have the same covariant derivative
in vertical directions, one has the associated covariant
derivatives coincide.\hfill$\Box$
\paragraph{}
Let us consider the pull-back bundle $\pi^* {\bf TM}\rightarrow
{\bf N}$ and the tangent bundle ${\bf TM\rightarrow {\bf M}}$
endowed with a linear connection $\nabla$. The horizontal lift of
$\nabla $ (or pull-back connection, ([7, pg 57])) is a connection
on ${\pi^* {\bf TM}}\rightarrow {\bf N}$ defined by the condition
\begin{equation}
(\pi^* \nabla)_{\iota (X)} \pi^* S =\pi^*(\nabla_X S),\quad
\tilde{X}\in {\bf TM}.
\end{equation}
The parameterized geodesics of both connections $\pi^* \nabla$ and
$\nabla$ are the same,
\begin{displaymath}
(\pi^* \nabla)_{\iota_u (X)} \pi^*_u X =0\quad
\Leftrightarrow\quad \nabla_X X=0,
\end{displaymath}
In order to prove that, let us choose a local coordinate system on ${\bf M}$ and let us
write the geodesics equations in local coordinates. To do that, we
need the connection coefficients of the above connections. In
particular, the possibly non-zero connection coefficients in the
natural coordinates induced from the local coordinate system are
such that
\begin{displaymath}
\nabla _{\partial_j} \partial_k =\Gamma ^i _{jk} \partial_i
\,\Rightarrow \, \pi^*\nabla _{\delta_j} \pi ^* \partial_k =\pi^*
(\Gamma ^i _{jk}
\partial_i)=(\Gamma ^i _{jk}\pi^*
\partial_i).
\end{displaymath}

\paragraph{}

Let us concentrate on Berwald spaces know. We can prove now the
following,
\begin{proposicion}
Let $\nabla^{ch}$ be the Chern connection of a Finsler structure
$({\bf M, F)}$, $\nabla^{b}$ the linear Berwald connection and $<\nabla^{ch}>$ its averaged connection. Then the
structure is Berwald iff $\pi^*<\nabla^{ch}>=\nabla^{ch}$.
\end{proposicion}
{\bf Proof 1:} If $\pi^* <\nabla^{ch}>=\nabla^{ch}$, since the
induced horizontal connection $\pi^* <\nabla^{ch}>$ has the same
coefficients that $<\nabla^{ch}>$ and they live on ${\bf M}$, the
structure $({\bf M},F)$ is Berwald.

Let us suppose that the structure is Berwald. Then
$\pi^*<\nabla^{ch}>=\pi^*\,<1>\,\nabla^{ch} = \nabla^{ch}$. This
relation is checked writing the action of the average covariant
derivative on arbitrary vector sections.

{\bf Proof 2:} An alternative proof of is the following. We know
that
\begin{displaymath}
Tor(\nabla^{ch} )=0\, \Rightarrow Tor(<\nabla^{ch}>)=0.
\end{displaymath}
 On the other hand, the
parameterized geodesics of $\pi^*<\nabla^{ch}>$ are the same than
the geodesics of $<\nabla^{ch}>$. But if the space is Berwald, the
geodesic equation of $<\nabla^{ch}>$ are the same than the
geodesic equation of $\nabla^{ch}$. From this fact it follows
$\pi^*<\nabla^{ch}>=\nabla^{ch}$, because both have zero torsion.
If $\pi^*<\nabla^{b}>=\nabla^{b}$, the Berwald connection lives on
{\bf M} and therefore the structure is Berwald. \hfill$\Box$

The following results is direct from Szab\'o's theorem,
\begin{proposicion}
Let $({\bf M}, F)$ be a Finsler structure. Then there is an affine
equivalent Riemannian structure $({\bf M},h)$ iff the structure is
Berwald.
\end{proposicion}
{\bf Proof:} if there is an affine equivalence Riemannian
structure $h$ such that its Levi-Civita connection $\nabla^h$ has
the same parameterized geodesics as the linear Berwald connection
$\nabla^b$ and both connection have also null torsion, then both
connections are the same ([3], section 5.4) and since the
connection coefficients $^h\Gamma^i _{ij}$ live in {\bf M}, the
structure is Berwald. Conversely, if $({\bf M},F)$ is Berwald, its
Berwald connection is metrizable ([6]).\hfill$\Box$.

Recall that for Berwald spaces $\nabla^b =\nabla^{ch}$. Then,
\begin{proposicion}
Let $({\bf M},F)$ be a Berwald structure. Then any Riemannian metric $h$
on {\bf M} such that $\nabla^b \pi^* h=0$ implies that the
associated Levi-Civita connection $\nabla^h$ leaves invariant the
indicatrix under horizontal parallel transport.
\end{proposicion}
{\bf Proof}: If the Riemannian structure $h$ is conserved by the
Berwald connection, $\nabla^b \pi^* h =0$.  This implies that
$<\nabla^b > h=0$. In addition, $<\nabla^b>$ is torsion free.
Therefore, $<\nabla^b>=\nabla^h$. If $\nabla^b$ leaves invariant
the indicatrix, also $\pi^* <\nabla^b>\,=\pi^* \nabla^h$ leaves
invariant the structure. \hfill$\Box$

There is a converse of this result,
\begin{proposicion}
Let $({\bf M},F)$ be a Finsler structure. Then if there is a
Riemannian metric $h$ that leaves invariant the indicatrix under
the parallel transport of
$\pi^*\nabla^h$, the structure is Berwald.
\end{proposicion}
{\bf Proof}: Let us consider such Riemannian metric $h$ and the
associated Levi-Civita connection $\nabla^h$. The induced
connection $ \pi^* \nabla ^h$ is torsion free, its connection
coefficients in natural coordinates live on {\bf M} and the
averaged connection $<\pi^*\nabla^h>$ coincides with $\nabla^h$,
so $ \pi^* \nabla ^h=\pi^* <\pi^*\nabla ^h>=\nabla^b $. The last
equality because $\pi^* <\pi^*\nabla ^h>$ leaves invariant the
indicatrix and it is torsion-free, therefore must be the Berwald
connection. Then the connection $\pi^* <\pi^*\nabla ^h>=\nabla^b$
has coefficients living on {\bf M} and the structure is
Berwald.\hfill$\Box$
\section{A corollary on non-Berwaldian Spaces Landsberg}

Let us consider a Riemannian metric $h$ such that its parallel Riemannian
transport leaves invariant the indicatrix of the Finsler metric
$F$, following {\it proposition 5.3.7}. Therefore $F$ is Berwald.
Let us also consider the set of interpolating metrics,
\begin{displaymath}
F_t (x,y)=(1-t)F(x,y)+t\sqrt{h(x)_{ij}y^i y^j},\,\,
i,j=1,...,n,\,\, t\in [0,1]
\end{displaymath}
and their indicatrix,
\begin{displaymath}
{\bf I}_x(t):=\{F_t (x,y)=1,\, y\in {\bf T}_x {\bf M},\,x \in {\bf M}\}.
\end{displaymath}
Since the metric $F$ is Berwald, each of the above interpolating
metrics defines an indicatrix which is invariant under the action
of the Levi-Civita connection of $h$: the parallel transport along
$\gamma(s)\subset {\bf M}$ of ${\bf I}_x$ leads to the indicatrix
over the final point of the path $\gamma$.

Let us that each of these indicatrix defines a submanifold of
${\bf T}_x{\bf M}$ of co-dimension 1 and that they are
non-intersecting sub-manifolds. Therefore the union of indicatrix
$\{{\bf I}_x (t),\,\in [0,1]\}$ defines a sub-manifold of ${\bf
T}_x{\bf M}$ of co-dimension $0$ that is invariant under the
holonomy of the metric $h$.
\paragraph{}
\begin{definicion}
A Finsler structure $({\bf M},F)$ is a Landsberg space if the
$hv$-curvature $P$ of the Chern's connection is such 
that $\dot{A}_{ijk}=P^n _{ijk}=0$, where the vector field is
defined as $e_n =\frac{y}{F(y)}$. A pure Landsberg space is such
that it is Landsberg and it is not Berwald.
\end{definicion}
This definition that we take of Landsberg space is a bit unusual,
although can be obtained from the standard characterizations
straightforwardly. In particular, the standard definition of
Landsberg space is such that ([1, {\it section 3.4}])
\begin{displaymath}
0=\dot{A}_{ikl}=-l^j\,P_{jikl}\,=\tilde{l}_j\,P^{j}_{ikl}:= P^{n}
_{ikl}.
\end{displaymath}
\begin{teorema}
Let $({\bf M},F)$ be a Finsler space and suppose that the averaged
connection $<\nabla^{ch}>$ does not leave invariant any compact
submanifolds ${\bf I}_x(t)\subset {\bf T}_x{\bf M}$ of codimension
zero. Then the structure $ ({\bf M},F)$ is a pure Landsberg space.
\end{teorema}
{\bf Proof:} suppose that the Landsberg space is Berwald. Then we
know from a theorem of Szabo that this linear Berwald connection
is metrizable ([6]). Then, there is a Riemannian connection
$\nabla^h$ that is identified with the average connection
$<\nabla^{ch}>$ and this is in contradiction with the hypothesis
of the theorem because $\pi^* \nabla^h =\pi^*
<\nabla^{ch}>=\nabla^{h}$ leaves invariant the set of indicatrix
${\bf I}_x(t),\, \forall t\in [0,1]$ as we show before, the union
defining a submanifold of co-dimension zero of ${\bf T}_x{\bf
M}$.\hfill$\Box$

\paragraph{}
One can use {\it theorem 5.4.2} to argue for a strategy to solve
the longest posed problem in Finsler Geometry. It is the {\it
conjecture} that there are not pure Landsberg spaces. The idea is
to use the classification of affine connections to show, using
additional techniques and constrains, that in fact, there are no
possible holonomies groups of affine connections ([4]) available.

The proof of the conjecture has been done by the author in
dimension $2$ using holonomy constrains but without using the
result of theorem {\it 5.4.2}. This is because in dimension $2$,
the number of possible averaged holonomies for Landsberg spaces is small and one
can check directly that in fact is not possible Landsberg spaces.
The problem is that in higher dimensions, the number of possible
holonomies grow. Therefore, the constrain of {\it theorem 5.4.2}
could play a role.

\chapter{Conclusions}

 From our point of view, the results presented in this work
reveals the power of a principle that we called {\it Convex
Invariance} in ref. [2]: the intrinsic invariance of some
geometric properties under a homotopy in the corresponding
operator {\it moduli space} of connections having the same averaged. Convex
 invariance is the invariance under a convex homotopy from any of these linear connections
to the average connection. The set of linear connections having the same
averaged is an equivalence class, which is a strongly convex set.
We say that a property is convex invariant if it is well define on each equivalence class.

We can see the results presented in this work from this perspective.
One example of how the principle works is the problem of
geodesic equivalence between different Finslerian structures. From
the point of view of Convex Invariance, one states the following
question:
\paragraph{}
{\it Which properties of the geodesics are defined on each
equivalence class?}
\paragraph{}
We prove that in the category of Berwald spaces, the
geodesics are the same on each equivalence class.

Another application of this point of view is the formulation of
the Landsberg problem in the following way:
\paragraph{}
{\it Is it the property of being Landsberg convex invariant?}
\paragraph{}
If the answer is yes, there are not Landsberg spaces which are not
Berwald.
\paragraph{}
Apart from the above results, the perspective adopted in the conclusion
seems applicable to other properties.

\end{document}